# DOUBLE INTEGRALS AND INFINITE PRODUCTS FOR SOME CLASSICAL CONSTANTS VIA ANALYTIC CONTINUATIONS OF LERCH'S TRANSCENDENT

JESÚS GUILLERA AND JONATHAN SONDOW

ABSTRACT. The two-fold aim of the paper is to unify and generalize on the one hand the double integrals of Beukers for $\zeta(2)$ and $\zeta(3)$, and of the second author for Euler's constant $\gamma$ and its alternating analog $\ln(4/\pi)$, and on the other hand the infinite products of the first author for $e$, of the second author for $\pi$, and of Ser for $e^\gamma$. We obtain new double integral and infinite product representations of many classical constants, as well as a generalization to Lerch's transcendent of Hadjicostas's double integral formula for the Riemann zeta function, and logarithmic series for the digamma and Euler beta functions. The main tools are analytic continuations of Lerch's function, including Hasse's series. We also use Ramanujan's polylogarithm formula for the sum of a particular series involving harmonic numbers, and his relations between certain dilogarithm values.

## Contents

## 1. INTRODUCTION

This paper is primarily about double integrals and infinite products related to the Lerch transcendent $\Phi(z, s, u)$.

Concerning double integrals, our aim is to unify and generalize the integrals over the unit square $[0, 1]^2$ of Beukers [4] and Hadjicostas [9], [10] for values of the Riemann zeta function (Examples 3.8 and 4.1), and those of the second author [20], [22] for Euler's constant $\gamma$ and its alternating analog $\ln(4/\pi)$ (Examples 4.4 and 3.14). We do this in Theorems 3.1 and 4.1, using a classical integral analytic continuation of $\Phi$ (Lemmas 2.1 and 2.2). As applications, we obtain new double integral representations of many constants, including $\ln 2$, $\ln \varphi$, $\pi^3$, $\ln \sigma$, $\Gamma(3/4)$, $\zeta(5)$, and $G/\pi$ (Examples 3.4, 3.5, 3.9, 3.13, 3.19, 3.21, and 3.24), where $\varphi$ is the golden ratio, $\sigma$ is one of Somos's quadratic recurrence constants, and $G$ is Catalan's constant. (Single integrals for them follow by the change of variables in the proof of Theorem 3.1.) Corollary 3.1 gives a double integral for the polylogarithm function; examples use Ramanujan's relations between certain values of the dilogarithm (Example 3.6). Corollaries 3.5, 3.6 and 4.1 evaluate certain types of double integrals in terms of the function $\ln \Gamma(x)$ and its derivative, the digamma function; Examples 3.14, 3.16 and 4.4 for $\ln \pi$, $\pi/\sqrt{3}$ and $\gamma$ are interesting special cases. Theorem 3.2 evaluates a class of double integrals with the help of Ramanujan's summation of a series involving harmonic numbers; applications use classical values of the di- and trilogarithm (Example 3.25).

Regarding infinite products, our aim is extend those of the first author (see [23]) for $e$ (Example 5.12), of the second author [23] for $\pi$ (Example 5.1), and of Ser [17] (see also [21], [23]) for $e^\gamma$ (Example 5.8). To do this we use two analytic continuations of $\Phi$ by series, one due to Hasse [12] (Theorems 2.1 and 2.2). We find new infinite products for many constants, among them $\Gamma(1/4)\Gamma(3/4)^{-1}$, $e^{G/\pi}$, $\pi/e$, $e^\pi$, $\sigma$, and the Glaisher-Kinkelin constant $A$ (Examples 5.4, 5.5, 5.6, 5.9, 5.10, and 5.11). The product for $\sigma$ has rational factors.

Finally, we mention other results obtained: explicit formulas for the Bernoulli and Euler polynomials (Examples 2.4 and 2.5), and logarithmic series for the digamma and Euler beta functions (Theorems 5.1 and 5.2). The latter series leads to infinite products for the exponential function of the same shape as other products in the paper, but which we construct without using Lerch's function (Theorem 5.3 and Examples 5.12 and 5.13).

## 2. The Lerch Transcendent

The *Lerch transcendent* $\Phi$ (see [2, section 1.11], [24, section 64:12]) is the analytic continuation of the series

$$(1) \qquad \Phi(z, s, u) = \frac{1}{u^s} + \frac{z}{(u+1)^s} + \frac{z^2}{(u+2)^s} + \cdots,$$

which converges for any real number $u > 0$ if $z$ and $s$ are any complex numbers with either $|z| < 1$, or $|z| = 1$ and $\Re(s) > 1$. Special cases include the analytic continuations of the *Riemann zeta function*

$$(2) \qquad \zeta(s) = \sum_{k=1}^{\infty} \frac{1}{k^s} = \Phi(1, s, 1),$$

the *Hurwitz zeta function*

$$\zeta(s, u) = \sum_{k=0}^{\infty} \frac{1}{(u+k)^s} = \Phi(1, s, u),$$

the *alternating zeta function* (also known as *Dirichlet's eta function* $\eta(s)$)

$$(3) \qquad \zeta^*(s) = \sum_{k=1}^{\infty} \frac{(-1)^{k-1}}{k^s} = \Phi(-1, s, 1),$$

the *Dirichlet beta function*

$$(4) \qquad \beta(s) = \sum_{k=0}^{\infty} \frac{(-1)^k}{(2k+1)^s} = 2^{-s}\Phi\left(-1, s, \frac{1}{2}\right),$$

the *Legendre chi function*

$$\chi_s(z) = \sum_{k=0}^{\infty} \frac{z^{2k+1}}{(2k+1)^s} = 2^{-s} z \Phi\left(z^2, s, \frac{1}{2}\right), \tag{5}$$

the *polylogarithm*

$$\text{Li}_n(z) = \sum_{k=1}^{\infty} \frac{z^k}{k^n} = z\Phi(z, n, 1), \tag{6}$$

and the *Lerch zeta function* $L(\lambda, \alpha, s) = \Phi(\exp(2\pi i \lambda), s, \alpha)$.

From (1), Lerch's transcendent satisfies the following identities:

$$\Phi(z, s, u+1) = \frac{1}{z}\left(\Phi(z, s, u) - \frac{1}{u^s}\right), \tag{7}$$

$$\Phi(z, s-1, u) = \left(u + z\frac{\partial}{\partial z}\right)\Phi(z, s, u), \tag{8}$$

$$\Phi(z, s+1, u) = -\frac{1}{s}\frac{\partial \Phi}{\partial u}(z, s, u). \tag{9}$$

Henceforth, we assume that $u$ and $v$ are positive real numbers. The following two lemmas are classical. (See [2], [15], [24], and [25].)

**Lemma 2.1.** *If either $|z| < 1$ and $\Re(s) > 0$, or $z = 1$ and $\Re(s) > 1$, then*

$$\Phi(z, s, u) = \frac{1}{\Gamma(s)} \int_0^{\infty} \frac{e^{-(u-1)t}}{e^t - z} t^{s-1} dt. \tag{10}$$

*Proof.* Replacing $\frac{e^{-(u-1)t}}{e^t - z}$ with $\sum_{k \geq 0} z^k e^{-(u+k)t}$, and noting that

$$\int_0^{\infty} e^{-(u+k)t} t^{s-1} dt = \frac{1}{(u+k)^s} \int_0^{\infty} e^{-x} x^{s-1} dx = \frac{\Gamma(s)}{(u+k)^s}, \tag{11}$$

the lemma follows. □

**Lemma 2.2.** *The series (1) extends to a function $\Phi(z, s, u)$ which is defined, and holomorphic in $z$ and $s$, for $z \in \mathbf{C} - [1, \infty)$ and all complex $s$, and which is given by the integral formula (10) if $\Re(s) > 0$.*

*Proof.* If $z \in \mathbf{C} - [1, \infty)$, the denominator of the integrand in (10) does not vanish. It follows that for $\Re(s) > 0$ the right side of (10) defines a function which is holomorphic in $z$ and $s$. The extension to other values of $s$ follows inductively using the identity (8). □

For $z \in \mathbf{C} - [1, \infty)$ and $s = 0, -1, -2, \ldots$, there is a closed formula for the function $\Phi(z, s, u)$. In fact, summing the series $\Phi(z, 0, u) = 1 + z + z^2 + \cdots$ and using (8), we see inductively that the functions

$$\Phi(z, 0, u) = \frac{1}{1-z}, \tag{12}$$

$$\Phi(z, -1, u) = \frac{u}{1-z} + \frac{z}{(1-z)^2}, \quad \cdots \tag{13}$$

are rational, with rational coefficients, holomorphic in $z$ except for a pole at $z = 1$. In particular, setting $z = -1$ and $m = -s$, we define the *mth Euler polynomial*

(14) $$E_m(x) = 2\Phi(-1, -m, x).$$

On the other hand, if $s$ and $u$ are positive integers, then $\Phi(z, s, u)$ is a transcendental function. For example, when $s = u = 1$ series (1) is $1 + \frac{z}{2} + \frac{z^2}{3} + \cdots$, which evidently sums to

(15) $$\Phi(z, 1, 1) = -\frac{1}{z}\ln(1-z) = \frac{1}{z}\text{Li}_1(z).$$

For $s = 2$, setting $u = 1$ in (8) and writing the right side as $\frac{\partial}{\partial z} z\Phi(z, 2, 1)$, we find that the integral

$$\Phi(z, 2, 1) = \frac{1}{z}\int_0^z \Phi(w, 1, 1)dw = \frac{1}{z}\text{Li}_2(z)$$

is $z^{-1}$ times the classical *dilogarithm function*. The identity (7) leads by induction to similar formulas for $\Phi(z, 1, u)$ and $\Phi(z, 2, u)$ when $u = 2, 3, \ldots$.

The following three examples give evaluations of $\partial\Phi/\partial s$ that we will need.

**Example 2.1.** Differentiating the relation (see, for example, [19])

(16) $$\zeta^*(s) = \left(1 - \frac{2}{2^s}\right)\zeta(s)$$

at $s = -1$, and using (3) and the values $\zeta(-1) = -1/12$ and

(17) $$\zeta'(-1) = \frac{1}{12} - \ln A,$$

where

$$A = \lim_{n\to\infty} \frac{1^1 2^2 3^3 \cdots n^n}{n^{\frac{1}{2}n^2+\frac{1}{2}n+\frac{1}{12}} e^{-\frac{1}{4}n^2}}$$

is the *Glaisher-Kinkelin constant* (see [7, section 2.15]), we get

(18) $$\frac{\partial\Phi}{\partial s}(-1, -1, 1) = \ln\frac{A^3}{2^{1/3}e^{1/4}}.$$

**Example 2.2.** Differentiating (16) and the functional equation of the zeta function [24, p. 27]

$$\zeta(s) = 2(2\pi)^{s-1}\Gamma(1-s)\sin\frac{\pi s}{2}\zeta(1-s)$$

at $s = -2$, and using (3) and the value $\zeta(-2) = 0$, we deduce that

(19) $$\frac{\partial\Phi}{\partial s}(-1, -2, 1) = \frac{7\zeta(3)}{4\pi^2},$$

where $\zeta(3)$ is *Apéry's constant*.

**Example 2.3.** Differentiating (4) and the functional equation of the Dirichlet beta function [24, p. 27],

$$\beta(s) = \left(\frac{\pi}{2}\right)^{s-1}\Gamma(1-s)\cos\frac{\pi s}{2}\beta(1-s)$$

at $s = -1$, and using the value $\beta(-1) = 0$, we obtain the evaluation

$$\frac{\partial \Phi}{\partial s}\left(-1, -1, \frac{1}{2}\right) = \frac{G}{\pi}, \tag{20}$$

where

$$G = \beta(2) = \frac{1}{4}\Phi\left(-1, 2, \frac{1}{2}\right) = \Re(\Phi(i, 2, 1)) \tag{21}$$

is *Catalan's constant*.

**Theorem 2.1.** *For all complex $s$, and complex $z$ with $\Re(z) < 1/2$,*

$$(1-z)\Phi(z, s, u) = \sum_{n=0}^{\infty} \left(\frac{-z}{1-z}\right)^n \sum_{k=0}^{n} (-1)^k \binom{n}{k} (u+k)^{-s}. \tag{22}$$

*Proof.* The right side of (22) defines a holomorphic function of $z$ on the half-plane $\Re(z) < 1/2$, because the inequality implies that $|\frac{-z}{1-z}| < 1$. We prove (22) when $|z| < 1/2$, and the result follows by analytic continuation.

Reversing the order of summation (for justification, see [6]), the right side of (22) is equal to

$$\sum_{k=0}^{\infty} \frac{(-1)^k}{(u+k)^s} \sum_{n=k}^{\infty} \binom{n}{k} \left(\frac{-z}{1-z}\right)^n = \sum_{k=0}^{\infty} \frac{(-1)^k}{(u+k)^s}(1-z)(-z)^k = (1-z)\Phi(z, s, u),$$

and the theorem follows. □

**Remark.** Theorem 2.1 excludes the Hurwitz zeta function $\zeta(s, u) = \Phi(1, s, u)$. However, the two cases $u = 1$ and $u = 1/2$ can be obtained via the relations (2), (3), (16), and $\zeta(s, 1/2) = (2^s - 1)\zeta(s, 1)$. In the case $u = 1$, formula (22) with $z = -1$ gives a globally convergent series for the alternating zeta function (3). Then (16) yields a formula for the Riemann zeta function, valid for all $s \neq 1$,

$$\zeta(s) = \frac{1}{1 - 2^{1-s}} \sum_{n=0}^{\infty} \frac{1}{2^{n+1}} \sum_{k=0}^{n} (-1)^k \binom{n}{k} (k+1)^{-s},$$

which was conjectured by Knopp, proved by Hasse [12, p. 464], and rediscovered by the second author [19].

**Theorem 2.2** (Hasse). *For $z = 1$, the series (1) extends to a function*

$$\Phi(1, s, u) = \frac{1}{s-1} \sum_{n=0}^{\infty} \frac{1}{n+1} \sum_{k=0}^{n} (-1)^k \binom{n}{k} (u+k)^{1-s} \tag{23}$$

*defined for all complex $s \neq 1$, with a simple pole at $s = 1$. Moreover, the series in (23) converges uniformly on compact sets in the $s$-plane to an entire function.*

*Proof.* We give a sketch; for details, see [12].

First assume that $\Re(s) > 1$. Setting $z = 1$ in (10), we rewrite it as

$$\Phi(1, s, u) = \frac{1}{\Gamma(s)} \int_0^{\infty} \frac{t}{1-e^{-t}} e^{-ut} t^{s-2} dt = \frac{1}{\Gamma(s)} \int_0^{\infty} \sum_{n=0}^{\infty} \frac{(1-e^{-t})^n}{n+1} e^{-ut} t^{s-2} dt.$$

Interchanging the integral and sum, we expand the binomial and get

$$\Phi(1, s, u) = \frac{1}{\Gamma(s)} \sum_{n=0}^{\infty} \frac{1}{n+1} \sum_{k=0}^{n} (-1)^k \binom{n}{k} \int_0^{\infty} e^{-(u+k)t} t^{s-2} dt.$$

Now (11) (with $k$ in place of $n$, and $s-1$ in place of $s$) and the relation $\Gamma(s) = (s-1)\Gamma(s-1)$ yield (23) when $\Re(s) > 1$.

For any $s \in \mathbf{C}$, the inner sum in (23) with $n > 0$ evidently is equal to

$$(24) \quad \sum_{k=0}^{n} (-1)^k \binom{n}{k} (u+k)^{1-s} = (s-1)_n \int_0^1 \cdots \int_0^1 (u + x_1 + \cdots x_n)^{1-s-n} dx_1 \cdots dx_n,$$

where the Pochhammer symbol $(a)_n$ stands for the product $a(a+1)\cdots(a+n-1)$. It follows (see [12] or [20]) that the series in (23) converges as required. By analytic continuation, this proves the theorem. □

**Example 2.4.** Note that when $s$ is zero or a negative integer, the factor $(s-1)_n$ in (24) vanishes for $n > m := 1 - s$, so that

$$(25) \quad \sum_{k=0}^{n} (-1)^k \binom{n}{k} (u+k)^m = 0 \qquad (n > m = 0, 1, 2, \dots).$$

In this case the series in (23) becomes a finite sum, in fact, a polynomial with rational coefficients. Defining the $m$th *Bernoulli polynomial* by $B_0(x) = 1$ and $B_m(x) = -m\Phi(1, 1-m, x)$ for $m = 1, 2, \dots$, we thus have the explicit formula

$$B_m(x) = \sum_{n=0}^{m} \frac{1}{n+1} \sum_{k=0}^{n} (-1)^k \binom{n}{k} (x+k)^m.$$

In particular,

$$(26) \quad B_1(x) = -\Phi(1, 0, x) = x - \frac{1}{2}.$$

**Example 2.5.** Taking $z = -1$ and $s = -m$ in (22), we use (25) to get the explicit formula for the $m$th Euler polynomial (14)

$$E_m(x) = \sum_{n=0}^{m} \frac{1}{2^n} \sum_{k=0}^{n} (-1)^k \binom{n}{k} (x+k)^m.$$

The next two results relate the analytic continuation of $\Phi(z, s, u)$ for $z \in \mathbf{C} - [1, \infty)$ (Lemmas 2.1 and 2.2) to that for $z = 1$ (Theorem 2.2). (Note that $\Phi(z, s, u)$ *is not continuous in $z$ at $z = 1$*: see (12), (13) and [2, p. 30, equation (12)].)

**Corollary 2.1.** *For all complex $s$,*

$$\int_0^{\infty} \frac{\Phi(-z, s, u)}{1+z} dz = s\Phi(1, s+1, u).$$

*Proof.* Replace $z$ with $-z$ in (22), multiply by $(1+z)^{-2}$, and integrate from $z = 0$ to $\infty$. Comparing the result with (23) gives the corollary. □

**Lemma 2.3.** *If either* $z \in \mathbf{C} - \{(-\infty, -1] \cup [1, \infty)\}$ *and* $s \in \mathbf{C}$, *or* $z = \pm 1$ *and* $s \in \mathbf{C} - \{1\}$, *then*

$$\Phi(z, s, u) = 2^{-s} \left[ \Phi\left(z^2, s, \frac{u}{2}\right) + z\Phi\left(z^2, s, \frac{u+1}{2}\right) \right]. \tag{27}$$

*Proof.* If $|z| < 1$ (respectively, $z = \pm 1$ and $\Re(s) > 1$), formula (27) is proved by splitting the series (1) into two sums, one over even numbers and the other over odd numbers, and factoring out $2^{-s}$. The result then follows by analytic continuation in $z$ (respectively, in $s$). □

As an application, (27) with $z = \pm 1$ and $s = 0, -1, -2, \ldots$ implies the classical relations among the Bernoulli and Euler polynomials [24, equations 19:5:6 and 20:3:3].

**Example 2.6.** Differentiating (27) with respect to $s$ at $s = 0$, and using (26) and the formula $\frac{\partial \Phi}{\partial s}(1, 0, u) = \ln \frac{\Gamma(u)}{\sqrt{2\pi}}$ (see [2, p. 26]), we get

$$\frac{\partial \Phi}{\partial s}(-1, 0, u) = \ln \frac{\Gamma\left(\frac{u}{2}\right)}{\Gamma\left(\frac{u+1}{2}\right)\sqrt{2}}. \tag{28}$$

## 3. Double Integrals

Using the integral analytic continuation of Lerch's function in Lemmas 2.1 and 2.2, we evaluate certain types of double integrals.

**Theorem 3.1.** *If either* $z \in \mathbf{C} - [1, \infty)$ *and* $\Re(s) > -2$, *or* $z = 1$ *and* $\Re(s) > -1$, *then*

$$\int_0^1 \int_0^1 \frac{x^{u-1}y^{v-1}}{1 - xyz}(-\ln xy)^s dx dy = \Gamma(s+1)\frac{\Phi(z, s+1, v) - \Phi(z, s+1, u)}{u - v}, \tag{29}$$

$$\int_0^1 \int_0^1 \frac{(xy)^{u-1}}{1 - xyz}(-\ln xy)^s dx dy = \Gamma(s+2)\Phi(z, s+2, u). \tag{30}$$

*Proof.* The integrals define holomorphic functions of $z$ and/or $s$ under the conditions stated. Making the change of variables $x = X/Y$, $y = Y$ in integral (29), it is equal to

$$\int_0^1 \frac{X^{u-1}}{1 - Xz}(-\ln X)^s \int_X^1 Y^{v-u-1} dY dX = \frac{1}{u - v} \int_0^1 \frac{X^{v-1} - X^{u-1}}{1 - Xz}(-\ln X)^s dX.$$

Substituting $t = -\ln X$ and using (10) (with $s + 1$ in place of $s$), we obtain equality (29) when $\Re(s) > 0$ and either $|z| < 1$ or $z = 1$; the result then follows by analytic continuation in $z$ or $s$. To prove (30), let $v \to u$ in (29), and use the identity (9). □

**Example 3.1.** Setting $z = -i$, $u = 1$ and $s = 0$ in (30), we use (1), (16), (21), and Euler's formula $\zeta(2) = \pi^2/6$ to get

$$\int_0^1 \int_0^1 \frac{1}{1 + xyi} dx dy = G - \frac{\pi^2 i}{48}.$$

The next three examples use the substitution $x \to x^2$, $y \to y^2$.

**Example 3.2.** Setting $z = -1$, $u = 1$, $v = 1/2$ and $s = 1$ in (29), we get

$$\int_0^1 \int_0^1 \frac{-x \ln xy}{1 + x^2 y^2} dx dy = G - \frac{\pi^2}{48}.$$

**Example 3.3.** Setting $z = 1$, $u = 1$, $v = 1/2$ and $s = 1$ in (29) gives

$$\int_0^1 \int_0^1 \frac{-x \ln xy}{1 - x^2 y^2} dx dy = \frac{\pi^2}{12}.$$

**Corollary 3.1.** *If either $z \in \mathbb{C} - [1, \infty)$ and $n \geq -1$, or $z = 1$ and $n \geq 0$, then*

$$\int_0^1 \int_0^1 \frac{(-\ln xy)^n}{1 - xyz} dx dy = \frac{(n+1)! \operatorname{Li}_{n+2}(z)}{z}.$$

*In particular, for $z \in \mathbb{C} - [1, \infty)$*

$$\int_0^1 \int_0^1 \frac{-1}{(1 - xyz) \ln xy} dx dy = \frac{-\ln(1-z)}{z}.$$

*Proof.* Let $u = 1$ in (30) and set $s = n$. Then use (6) and (15). □

**Example 3.4.** Take $z = 1/2$. Using Euler's and Landen's values of the di- and trilogarithm at $1/2$ (see [14, pp. 1-2]), we see that (compare Example 3.13)

$$\int_0^1 \int_0^1 \frac{-1}{(2 - xy) \ln xy} dx dy = \ln 2, \qquad \int_0^1 \int_0^1 \frac{1}{2 - xy} dx dy = \frac{\pi^2}{12} - \frac{\ln^2 2}{2},$$

$$\int_0^1 \int_0^1 \frac{-\ln xy}{2 - xy} dx dy = \frac{7\zeta(3)}{4} - \frac{\pi^2 \ln 2}{6} + \frac{\ln^3 2}{3}.$$

**Example 3.5.** Take $z = \varphi^{-1}$, $\varphi^{-2}$, $-\varphi$, and $-\varphi^{-1}$, where $\varphi = (1 + \sqrt{5})/2$ is the *golden ratio*. Using the relation $1 + \varphi = \varphi^2$ to simplify $\ln(1 - z)$, and using Landen's values for $\operatorname{Li}_2(z)$ and $\operatorname{Li}_3(\varphi^{-2})$ (see [13, equations (1.20), (1.21), and (6.13)]), we obtain

$$\int_0^1 \int_0^1 \frac{-1}{(\varphi - xy) \ln xy} dx dy = 2 \ln \varphi, \qquad \int_0^1 \int_0^1 \frac{1}{\varphi - xy} dx dy = \frac{\pi^2}{10} - \ln^2 \varphi,$$

$$\int_0^1 \int_0^1 \frac{-1}{(\varphi^2 - xy) \ln xy} dx dy = \ln \varphi, \qquad \int_0^1 \int_0^1 \frac{1}{\varphi^2 - xy} dx dy = \frac{\pi^2}{15} - \ln^2 \varphi,$$

$$\int_0^1 \int_0^1 \frac{-1}{(1 + \varphi xy) \ln xy} dx dy = \frac{2 \ln \varphi}{\varphi}, \qquad \int_0^1 \int_0^1 \frac{1}{1 + \varphi xy} dx dy = \frac{\pi^2}{10\varphi} + \frac{\ln^2 \varphi}{\varphi},$$

$$\int_0^1 \int_0^1 \frac{-1}{(\varphi + xy) \ln xy} dx dy = \ln \varphi, \qquad \int_0^1 \int_0^1 \frac{1}{\varphi + xy} dx dy = \frac{\pi^2}{15} - \frac{\ln^2 \varphi}{2},$$

$$\int_0^1 \int_0^1 \frac{-\ln xy}{\varphi^2 - xy} dx dy = \frac{8\zeta(3)}{5} - \frac{4\pi^2 \ln \varphi}{15} + \frac{4 \ln^3 \varphi}{3}.$$

**Example 3.6.** Using Ramanujan's values [3, Part IV, p. 324] for $\operatorname{Li}_2(-1/8) + \operatorname{Li}_2(1/9)$, $\operatorname{Li}_2(-1/2) + \frac{1}{6} \operatorname{Li}_2(1/9)$, and $\operatorname{Li}_2(1/3) - \frac{1}{6} \operatorname{Li}_2(1/9)$, and substituting $x, y \to x^2, y^2$ in our double integral for the last $\text{Li}_2(1/9)$, we get, respectively,

$$\int_0^1 \int_0^1 \frac{1-2xy}{(8+xy)(9-xy)} dxdy = \frac{1}{2} \ln^2 \frac{9}{8},$$

$$\int_0^1 \int_0^1 \frac{52-7xy}{(2+xy)(9-xy)} dxdy = \frac{\pi^2}{3} + 3\ln^2 2 + 2\ln^2 3 - 6\ln 2 \ln 3,$$

$$\int_0^1 \int_0^1 \frac{9+xy}{9-x^2y^2} dxdy = \frac{\pi^2}{6} - \frac{\ln^2 3}{2}.$$

**Example 3.7.** In Bailey, Borwein, and Plouffe's [1] dilogarithm "ladder" $6\text{Li}_2(1/2) - 6\text{Li}_2(1/4) - 2\text{Li}_2(1/8) + \text{Li}_2(1/64) = \zeta(2)$, we replace $\text{Li}_2(1/2^k)$ with a double integral in which we substitute $x, y \to x^k, y^k$, for $k = 1, 2, 3$, and 6, obtaining

$$\int_0^1 \int_0^1 \frac{(4+2xy+x^2y^2)(4-8xy+x^2y^2)}{64-x^6y^6} dxdy = \frac{\pi^2}{72}.$$

**Corollary 3.2.** *We have*

(31) $$\int_0^1 \int_0^1 \frac{(-\ln xy)^s}{1-xy} dxdy = \Gamma(s+2)\zeta(s+2) \qquad (\Re(s) > -1),$$

(32) $$\int_0^1 \int_0^1 \frac{(-\ln xy)^s}{1+xy} dxdy = \Gamma(s+2)\zeta^*(s+2) \qquad (\Re(s) > -2),$$

(33) $$\int_0^1 \int_0^1 \frac{(-\ln xy)^s}{1+x^2y^2} dxdy = \Gamma(s+2)\beta(s+2) \qquad (\Re(s) > -2).$$

*Proof.* In (30) take $z = 1$, $u = 1$; $z = -1$, $u = 1$; and $z = -1$, $u = 1/2$; then use relations (2), (3) and (4), respectively. After using (4), substitute $x, y \to x^2, y^2$. □

For integer $s \geq 0$, formula (31) is due to Hadjicostas [9].

**Example 3.8.** Taking $s = 0$ and 1 in (31), we recover Beukers's integrals [4]

(34) $$\int_0^1 \int_0^1 \frac{1}{1-xy} dxdy = \zeta(2), \qquad \int_0^1 \int_0^1 \frac{-\ln xy}{1-xy} dxdy = 2\zeta(3).$$

**Example 3.9.** Taking $s = -1, 0$, and 1 in (33), and using (21) and the values $\beta(1) = \pi/4$ and $\beta(3) = \pi^3/32$ (see [24, p. 29]), we get

$$\int_0^1 \int_0^1 \frac{-1}{(1+x^2y^2)\ln xy} dxdy = \frac{\pi}{4}, \qquad \int_0^1 \int_0^1 \frac{1}{1+x^2y^2} dxdy = G,$$

$$\int_0^1 \int_0^1 \frac{-\ln xy}{1+x^2y^2} dxdy = \frac{\pi^3}{16}.$$

**Example 3.10.** Let $\rho$ be a non-trivial zero of the zeta function (2). Differentiate (32) at $s = \rho - 2$ and use (16) to get (compare Example 4.3)

$$\int_0^1 \int_0^1 \frac{(-\ln xy)^{\rho-2}}{1+xy} \ln(-\ln xy) dxdy = (1 - 2^{3-\rho})\Gamma(\rho)\zeta'(\rho).$$

**Corollary 3.3.** *If either $z \in \mathbf{C} - ((-\infty, -1] \bigcup [1, \infty))$ and $\Re(s) > -2$, or $z \pm 1$ and $\Re(s) > -1$, then*
$$\int_0^1 \int_0^1 \frac{(-\ln xy)^s}{1 - x^2y^2z^2} dx dy = \Gamma(s+2)\frac{\chi_{s+2}(z)}{z},$$
*where $\chi_s(z)$ is Legendre's chi function (5).*

*Proof.* In (30), take $u = 1/2$ and substitute $x, y, z \to x^2, y^2, z^2$. Then use (5). □

**Example 3.11.** Set $s = 0$, and take $z = \tan\frac{\pi}{8} = \sqrt{2} - 1$ and $z = \varphi^{-3} = \sqrt{5} - 2$. Using Landen's values [13, section 1.8.2] for $\chi_2(\tan\frac{\pi}{8})$ and $\chi_2(\varphi^{-3})$, we obtain, respectively,
$$\int_0^1 \int_0^1 \frac{1}{1 - x^2y^2\tan^2\frac{\pi}{8}} dx dy = \frac{\pi^2}{16\tan\frac{\pi}{8}} - \frac{\ln^2 \tan\frac{\pi}{8}}{4\tan\frac{\pi}{8}},$$
$$\int_0^1 \int_0^1 \frac{1}{\varphi^6 - x^2y^2} dx dy = \frac{\pi^2}{24\varphi^3} - \frac{3\ln^2\varphi}{4\varphi^3}.$$

**Corollary 3.4.** *For $z \in \mathbf{C} - [1, \infty)$,*

(35) $$\int_0^1 \int_0^1 \frac{-x^{u-1}y^{v-1}}{(1 - xyz)\ln xy} dx dy = \frac{1}{u-v}\left[\frac{\partial \Phi}{\partial s}(z, 0, v) - \frac{\partial \Phi}{\partial s}(z, 0, u)\right],$$

(36) $$\int_0^1 \int_0^1 \frac{-(xy)^{u-1}}{(1 - xyz)\ln xy} dx dy = \Phi(z, 1, u).$$

*Proof.* Substitute $(s+1)^{-1}\Gamma(s+2)$ for $\Gamma(s+1)$ in (29) and let $s \to -1^+$. Since $z \in \mathbf{C} - [1, \infty)$, in view of (12) the result is (35). Letting $v \to u$ in (35), and using identity (9), we arrive at (36). □

**Example 3.12.** Set $z = 0$. Using (1), we get
$$\int_0^1 \int_0^1 \frac{x^{u-1}y^{v-1}}{-\ln xy} dx dy = \frac{1}{u-v}\ln\frac{u}{v}, \qquad \int_0^1 \int_0^1 \frac{(xy)^{u-1}}{-\ln xy} dx dy = \frac{1}{u}.$$

**Example 3.13.** Set $z = 1/2$, $u = 2$, and $v = 1$ in (35). Using the derivative of (7) with respect to $s$ at $s = 0$, we get (compare Example 3.4)
$$\int_0^1 \int_0^1 \frac{-x}{(2 - xy)\ln xy} dx dy = \ln \sigma,$$
where

(37) $$\sigma = \sqrt{1\sqrt{2\sqrt{3\cdots}}} = 1^{1/2}2^{1/4}3^{1/8}\cdots = \exp\left[-\frac{1}{2}\frac{\partial \Phi}{\partial s}\left(\frac{1}{2}, 0, 1\right)\right]$$

is one of *Somos's quadratic recurrence constants* [18] (see [7, p. 446], which uses the notation $\gamma$ instead of $\sigma$; see also [16, p. 348]).

**Corollary 3.5.** *For $u > 0$ and $v > 0$,*

(38) $$\int_0^1 \int_0^1 \frac{-x^{u-1}y^{v-1}}{(1 + xy)\ln xy} dx dy = \frac{1}{u-v}\ln\frac{\Gamma(\frac{v}{2})\,\Gamma(\frac{u+1}{2})}{\Gamma(\frac{u}{2})\,\Gamma(\frac{v+1}{2})},$$

(39) $$\int_0^1 \int_0^1 \frac{-(xy)^{u-1}}{(1 + xy)\ln xy} dx dy = \frac{1}{2}\left[\psi\left(\frac{u+1}{2}\right) - \psi\left(\frac{u}{2}\right)\right],$$

where $\psi(u) = \frac{d}{du} \ln \Gamma(u) = \Gamma'(u)/\Gamma(u)$ *is the* digamma function.

*Proof.* To prove (38), set $z = -1$ in (35) and use (28). To prove (39), let $v \to u$ in (38). □

**Example 3.14.** Setting $u = 1$ in (39) and $u = 2$ and $v = 1$ in (38) gives integrals for $\ln 2$ and $\ln(\pi/2)$, respectively. Their sum and difference are integrals for $\ln \pi$ and the "alternating Euler constant" $\ln(4/\pi)$ (compare Example 4.4 and see [22] and [23]):

$$\int_0^1 \int_0^1 \frac{1+x}{(1+xy)(-\ln xy)} dxdy = \ln \pi, \qquad \int_0^1 \int_0^1 \frac{1-x}{(1+xy)(-\ln xy)} dxdy = \ln \frac{4}{\pi}.$$

**Example 3.15.** Set $u = 1$ and $v = 1/2$ in (38), and substitute $x \to x^2$, $y \to y^2$. Using the reflection formula $\Gamma(t)\Gamma(1-t) = \pi/\sin \pi t$ with $t = 1/4$, we get

$$\int_0^1 \int_0^1 \frac{-x}{(1+x^2y^2)\ln xy} dxdy = \ln \frac{\sqrt{2\pi}}{\Gamma(3/4)^2}.$$

**Corollary 3.6.** *For* $u > 0$ *and* $v > 0$,

(40) $$\int_0^1 \int_0^1 \frac{x^{u-1}y^{v-1}}{1-xy} dxdy = \frac{\psi(u) - \psi(v)}{u - v},$$

(41) $$\int_0^1 \int_0^1 \frac{(xy)^{u-1}}{1-xy} dxdy = \psi'(u).$$

*Proof.* To prove (40), set $z = 1$ in (29), let $s \to 0$, and use the formula [25, p. 271]

(42) $$\lim_{s \to 1} \left( \Phi(1, s, u) - \frac{1}{s-1} \right) = -\psi(u).$$

To prove (41), let $v \to u$ in (40). □

**Example 3.16.** Take $u = 1/3$ and $v = 2/3$ in (40). Since $\psi(2/3) - \psi(1/3) = \pi/\sqrt{3}$ (see [24, equation 44:7:1]), the change of variables $x \to x^3$ and $y \to y^3$ yields

$$\int_0^1 \int_0^1 \frac{y}{1-x^3y^3} dxdy = \frac{\pi}{3\sqrt{3}}.$$

**Example 3.17.** Taking $u = 1/2$ in (41) and using the value $\psi'(1/2) = \pi^2/2$ (see [24, p. 427]), we substitute $x \to x^2$, $y \to y^2$ and obtain

$$\int_0^1 \int_0^1 \frac{1}{1-x^2y^2} dxdy = \frac{\pi^2}{8}.$$

**Corollary 3.7.** *If* $\Re(s) > -2$, *then*

(43) $$\int_0^1 \int_0^1 x^{u-1}y^{v-1}(-\ln xy)^s dxdy = \Gamma(s+1)\frac{v^{-s-1} - u^{-s-1}}{u - v},$$

(44) $$\int_0^1 \int_0^1 (xy)^{u-1}(-\ln xy)^s dxdy = \Gamma(s+2)u^{-s-2}.$$

*Proof.* Set $z = 0$ in (29) and (30), and use (1). □

**Example 3.18.** Taking $s = -3/2$, $u = 2$ and $v = 1$ in (43) gives
$$\int_0^1 \int_0^1 \frac{x}{(-\ln xy)^{3/2}} dxdy = 2(\sqrt{2}-1)\sqrt{\pi}.$$

**Example 3.19.** Taking $u = 1$ and $s = -3/2$ or $s = -5/4$ in (44) gives
$$\int_0^1 \int_0^1 \frac{1}{(-\ln xy)^{3/2}} dxdy = \sqrt{\pi}, \qquad \int_0^1 \int_0^1 \frac{1}{(-\ln xy)^{5/4}} dxdy = \Gamma(3/4).$$

**Corollary 3.8.** *If $z \in \mathbf{C} - [1, \infty)$ and $\Re(s) > -2$, then*

$$(45) \quad \int_0^1 \int_0^1 \frac{x^{u-1} y^{v-1}}{(1-xyz)^2} (-\ln xy)^s dxdy$$
$$= \frac{\Gamma(s+1)}{u-v}[(1-v)\Phi(z, s+1, v) + \Phi(z, s, v) + (u-1)\Phi(z, s+1, u) - \Phi(z, s, u)],$$

$$(46) \quad \int_0^1 \int_0^1 \frac{(xy)^{u-1}}{(1-xyz)^2} (-\ln xy)^s dxdy = \Gamma(s+2)[(1-u)\Phi(z, s+2, u) + \Phi(z, s+1, u)].$$

*Proof.* Differentiate (29) with respect to $z$ and apply identity (8). Replacing $u$ with $u-1$, and $v$ with $v-1$, we arrive at (45). To prove (46), let $v \to u$ in (45) and use identity (9). □

**Example 3.20.** Set $z = -1$, $u = 1/2$, $v = 1$, and $s = 2$ in (45). Substituting $x \to x^2$, $y \to y^2$, and using (3), (4), (21), and the values $\zeta^*(2) = \pi^2/12$ and $\beta(3) = \pi^3/32$, we get
$$\int_0^1 \int_0^1 \frac{x \ln^2 xy}{(1+x^2y^2)^2} dxdy = G - \frac{\pi^2}{48} + \frac{\pi^3}{32}.$$

**Example 3.21.** Set $z = -1$, $u = 1$ and $s = 4$ in (46) to get
$$\int_0^1 \int_0^1 \frac{\ln^4 xy}{(1+xy)^2} dxdy = \frac{225}{2} \zeta(5).$$

**Example 3.22.** Setting $z = -1$, $u = 1/2$ and $s = -1$ in (46), and substituting $x \to x^2$, $y \to y^2$, we get
$$\int_0^1 \int_0^1 \frac{-1}{(1+x^2y^2)^2 \ln xy} dxdy = \frac{\pi+2}{8}.$$

**Corollary 3.9.** *If $z \in \mathbf{C} - [1, \infty)$, then*

$$\int_0^1 \int_0^1 \frac{-x^{u-1} y^{v-1}}{(1-xyz)^2 \ln xy} dxdy$$
$$= \frac{1}{u-v} \left[(1-v)\frac{\partial \Phi}{\partial s}(z, 0, v) + \frac{\partial \Phi}{\partial s}(z, -1, v) - (1-u)\frac{\partial \Phi}{\partial s}(z, 0, u) - \frac{\partial \Phi}{\partial s}(z, -1, u)\right].$$

*Proof.* Substitute $(s+1)^{-1}\Gamma(s+2)$ for $\Gamma(s+1)$ in (45) and let $s \to -1^+$. □

**Example 3.23.** Setting $z = -1$, $u = 2$ and $v = 1$, we use (28), (18), and (7) to get (compare [22, equation (9)])

$$\int_0^1 \int_0^1 \frac{-x}{(1+xy)^2 \ln xy} dx dy = \ln \frac{A^6}{2^{1/6}\sqrt{\pi e}},$$

where $A$ is the Glaisher-Kinkelin constant.

**Example 3.24.** Set $z = -1$, $u = 2$ and $v = 1/2$. Using (20) and substituting $x \to x^2$ and $y \to y^2$, we get

$$\int_0^1 \int_0^1 \frac{-x^2}{(1+x^2y^2)^2 \ln xy} dx dy = \frac{G}{\pi}.$$

For the proof of the final theorem in this section, we require a lemma.

**Lemma 3.1.** *If $H_{n,r} = 1 + \frac{1}{2^r} + \cdots + \frac{1}{n^r}$ is the nth harmonic number of order $r$ and $H_n = H_{n,1}$, and if $|z| < 1$, then*

$$(47) \quad \sum_{n=1}^{\infty} \frac{H_{n,2}}{n} z^n + 2 \sum_{n=1}^{\infty} \frac{H_n}{n^2} z^n = 3\text{Li}_3(z) - \text{Li}_2(z) \ln(1-z).$$

*Proof.* We have

$$\sum_{n=1}^{\infty} \frac{H_n}{n^2} z^n - \text{Li}_3(z) = \sum_{n=2}^{\infty} \sum_{k=1}^{n-1} \frac{z^n}{n^2 k} = \sum_{k=1}^{\infty} \sum_{n=k+1}^{\infty} \frac{z^n}{n^2 k} = \sum_{k=1}^{\infty} \sum_{j=1}^{\infty} \frac{z^{k+j}}{(k+j)^2 k}$$

$$= \sum_{k=1}^{\infty} \frac{1}{k} \sum_{j=1}^{\infty} \int_0^z \frac{1}{t} \int_0^t u^{k+j-1} du dt = \sum_{k=1}^{\infty} \frac{1}{k} \int_0^z \frac{1}{t} \int_0^t \frac{u^k}{1-u} du dt$$

$$= \int_0^z \frac{1}{t} \int_0^t \frac{-\ln(1-u)}{1-u} du dt = \frac{1}{2} \int_0^z \frac{\ln^2(1-t)}{t} dt$$

and in a similar way

$$\sum_{n=1}^{\infty} \frac{H_{n,2}}{n} z^n - \text{Li}_3(z) = \int_0^z \frac{\text{Li}_2(t)}{1-t} dt.$$

Adding twice the first identity to the second, we obtain

$$\sum_{n=1}^{\infty} \frac{H_{n,2}}{n} z^n + 2 \sum_{n=1}^{\infty} \frac{H_n}{n^2} z^n - 3\text{Li}_3(z) = \int_0^z \left[ \frac{\text{Li}_2(t)}{1-t} + \frac{\ln^2(1-t)}{t} \right] dt = -\text{Li}_2(z) \ln(1-z),$$

the evaluation of the integral being easily checked by differentiation. □

**Theorem 3.2.** *For $z \in \mathbf{C} - [1, \infty)$, the following analytic continuations hold:*

$$(48) \quad \int_0^1 \int_0^1 \frac{-\ln(1-xz)}{1-xy} dx dy$$
$$= \left[ \frac{1}{2} \ln z \ln(1-z) + \text{Li}_2(1-z) \right] \ln(1-z) + \text{Li}_3(z) - \text{Li}_3(1-z) + \zeta(3),$$

$$(49) \quad \int_0^1 \int_0^1 \frac{-\ln(1-xyz)}{1-xy} dx dy = \text{Li}_2(1-z) \ln(1-z) - \text{Li}_3(z) - 2\text{Li}_3(1-z) + 2\zeta(3).$$

*Proof.* By analytic continuation, it suffices to prove (48) and (49) when $|z| < 1$.

Expanding $-\ln(1-xz)$ in a series, and using (40) and the identity $\psi(n+1)-\psi(1) = H_n$ (see [24, equation 44:5:5]), we see that

$$\int_0^1 \int_0^1 \frac{-\ln(1-xz)}{1-xy} dxdy = \sum_{n=1}^\infty \frac{z^n}{n} \int_0^1 \int_0^1 \frac{x^n}{1-xy} dxdy = \sum_{n=1}^\infty \frac{H_n}{n^2} z^n.$$

Ramanujan's summation of the last series [3, p. 251, Entry 9(i) and p. 259, equation (12.1)] completes the proof of (48).

Similarly, using (41) and the identity $\psi'(n+1) = \frac{1}{6}\pi^2 - H_{n,2}$ (see [24, equation 44:12:6]), we find that

$$\int_0^1 \int_0^1 \frac{-\ln(1-xyz)}{1-xy} dxdy = -\frac{\pi^2}{6} \ln(1-z) - \sum_{n=1}^\infty \frac{H_{n,2}}{n} z^n.$$

Now apply (47) with $\sum_{n=1}^\infty H_n n^{-2} z^n$ replaced by the right side of (48). Simplifying with Euler's identity $\mathrm{Li}_2(z) + \mathrm{Li}_2(1-z) = \frac{1}{6}\pi^2 - \ln z \ln(1-z)$ (see [14, p. 1]), we arrive at (49). □

**Example 3.25.** First take $z = 1/2$; write $\ln\left(1 - \frac{x}{2}\right)$ as $\ln(2-x) - \ln 2$ in (48), rewrite (49) similarly, and apply the first equation in (34). Then take $z = -1$. Using the values of $\mathrm{Li}_2(1/2)$ and $\mathrm{Li}_3(1/2)$ (see Example 3.5), and the inversion formulas for $\mathrm{Li}_2(y) + \mathrm{Li}_2(1/y)$ and $\mathrm{Li}_3(y) - \mathrm{Li}_3(1/y)$ (see [13, equations (1.10) and (6.7)]) with $y = 2$, we obtain

$$\int_0^1 \int_0^1 \frac{\ln(2-x)}{1-xy} dxdy = \int_0^1 \int_0^1 \frac{\ln(1+xy)}{1-xy} dxdy = \frac{\pi^2 \ln 2}{4} - \zeta(3),$$

$$\int_0^1 \int_0^1 \frac{\ln(2-xy)}{1-xy} dxdy = \int_0^1 \int_0^1 \frac{\ln(1+x)}{1-xy} dxdy = \frac{5}{8}\zeta(3).$$

**Example 3.26.** Letting $z \to 1^-$, we deduce that

$$\int_0^1 \int_0^1 \frac{-\ln(1-x)}{1-xy} dxdy = 2 \int_0^1 \int_0^1 \frac{-\ln(1-xy)}{1-xy} dxdy = 2\zeta(3).$$

## 4. A GENERALIZATION OF HADJICOSTAS'S FORMULA

Hadjicostas [11] asked for a generalization of the double integral formulas in [10] and [22] for the Riemann and alternating zeta functions (Examples 4.1 and 4.2). To fill the bill, we multiply the integrand of (30) by $1 - x$, and obtain a formula valid on a larger half-plane of the complex $s$-plane.

**Theorem 4.1.** *If either $z \in \mathbf{C} - [1, \infty)$ and $\Re(s) > -3$, or $z = 1$ and $\Re(s) > -2$, then*

(50) $$\int_0^1 \int_0^1 \frac{1-x}{1-xyz} (xy)^{u-1}(-\ln xy)^s dxdy$$
$$= \Gamma(s+2) \left[ \Phi(z, s+2, u) + \frac{(1-z)\Phi(z, s+1, u) - u^{-s-1}}{z(s+1)} \right].$$

*Proof.* The integral in (50) defines a function which is holomorphic in $s$, when $\Re(s) > -3$ if $z \in \mathbf{C} - [1, \infty)$, and when $\Re(s) > -2$ if $z = 1$. We prove (50) with $\Re(s) > 0$ and the result then follows by analytic continuation.

Replace $u$ with $u+1$ in (29) and set $v = u$. Subtracting the result from (30) gives

$$(51) \quad \int_0^1 \int_0^1 \frac{(xy)^{u-1} - x^u y^{u-1}}{1 - xyz}(-\ln xy)^s dx dy$$
$$= \Gamma(s+2)\Phi(z, s+2, u) + \Gamma(s+1)\left[\Phi(z, s+1, u+1) - \Phi(z, s+1, u)\right].$$

Using identity (7), we obtain (50), and the theorem follows. □

**Example 4.1.** The case $z = u = 1$ is Hadjicostas's formula [5], [10] (see also [22])

$$(52) \quad \int_0^1 \int_0^1 \frac{1-x}{1-xy}(-\ln xy)^s dx dy = \Gamma(s+2)\left[\zeta(s+2) - \frac{1}{s+1}\right] \quad (\Re(s) > -2).$$

**Example 4.2.** Taking $z = -1$ and $u = 1$ gives the analogous formula [22]

$$\int_0^1 \int_0^1 \frac{1-x}{1+xy}(-\ln xy)^s dx dy = \Gamma(s+2)\left[\zeta^*(s+2) + \frac{1 - 2\zeta^*(s+1)}{s+1}\right] \quad (\Re(s) > -3).$$

**Example 4.3.** Let $\rho$ be a non-trivial zero of $\zeta(s)$. Differentiating (52) with respect to $s$ at $s = \rho - 2$, we get (compare Example 3.10)

$$\int_0^1 \int_0^1 \frac{1-x}{1-xy}(-\ln xy)^{\rho-2}\ln(-\ln xy) dx dy = \Gamma(\rho)\left[\zeta'(\rho) + \frac{1}{(\rho-1)^2} - \frac{\psi(\rho)}{\rho-1}\right].$$

**Corollary 4.1.** *For $u > 0$,*

$$\int_0^1 \int_0^1 \frac{1-x}{(1-xy)(-\ln xy)}(xy)^{u-1} dx dy = \ln u - \psi(u).$$

*Proof.* Set $z = 1$ in (50), let $s \to -1^+$, and use (42). □

**Example 4.4.** Setting $u = 1$ and using the value [24, p. 427]

$$(53) \quad \psi(1) = -\gamma,$$

where $\gamma = \lim_{n \to \infty}\left(1 + \frac{1}{2} + \cdots + \frac{1}{n} - \ln n\right)$ is *Euler's constant*, we recover the formula in [20], [22], and [23]

$$\int_0^1 \int_0^1 \frac{1-x}{(1-xy)(-\ln xy)} dx dy = \gamma.$$

**Corollary 4.2.** *For $z \in \mathbf{C} - [1, \infty)$,*

$$\int_0^1 \int_0^1 \frac{1-x}{(1-xyz)\ln^2 xy}(xy)^{u-1} dx dy$$
$$= \frac{\partial \Phi}{\partial s}(z, 0, u) + \frac{1}{z}\left[(z-1)\frac{\partial \Phi}{\partial s}(z, -1, u) - u\ln u\right] + \frac{1}{z-1}.$$

*Proof.* Let $s \to -2^+$ in (50) and use (12). □

**Example 4.5.** Take $u = 1/2$ and $z = -1$. Substituting $x \to x^2$, $y \to y^2$ and using (28) and (20), we get

$$\int_0^1 \int_0^1 \frac{1-x^2}{(1+x^2y^2)\ln^2 xy} dx dy = \ln \frac{\Gamma(1/4)}{2\Gamma(3/4)} + \frac{2G}{\pi} - \frac{1}{2}.$$

**Example 4.6.** Replacing $\Phi$ with the series (1), we let $z \to 0$ and deduce that

$$\int_0^1 \int_0^1 \frac{1-x}{\ln^2 xy} (xy)^{u-1} dx dy = (u+1)\ln\left(1+\frac{1}{u}\right) - 1.$$

**Corollary 4.3.** *If $\Re(s) > -3$, then*

$$\int_0^1 \int_0^1 (1-x)(xy)^{u-1}(-\ln xy)^s dx dy = \Gamma(s+1)\left[(s+1-u)u^{-s-2} + (u+1)^{-s-1}\right].$$

*Proof.* Set $z = 0$ in (51) and use (1). □

**Example 4.7.** Taking $s = -5/2$ and $u = 1$ gives

$$\int_0^1 \int_0^1 \frac{1-x}{(-\ln xy)^{5/2}} dx dy = \frac{\sqrt{\pi}}{3}(8\sqrt{2} - 10).$$

## 5. Infinite Products

Using corollaries of Theorems 2.1 and 2.2, together with logarithmic series for the digamma and Euler beta functions (Theorems 5.1 and 5.2), we derive infinite products for many constants.

**Corollary 5.1.** *For $m = 0, 1, 2, \ldots$ and complex $z$ with $\Re(z) < 1/2$,*

$$\frac{\partial \Phi}{\partial s}(z, -m, u) = \left(u + z\frac{\partial}{\partial z}\right)^m \sum_{n=0}^{\infty} \frac{1}{1-z}\left(\frac{-z}{1-z}\right)^n \sum_{k=0}^n (-1)^{k+1} \binom{n}{k} \ln(u+k).$$

*Proof.* Differentiate (22) with respect to $s$ at $s = 0$ and multiply by $(1-z)^{-1}$. Then apply (8) $m$ times. □

**Example 5.1.** Setting $m = 0$, $z = -1$ and $u = 1$ and multiplying by 2 gives

$$2\frac{\partial \Phi}{\partial s}(-1, 0, 1) = \sum_{n=0}^{\infty} \frac{1}{2^n} \sum_{k=0}^n (-1)^{k+1} \binom{n}{k} \ln(k+1).$$

Using (28), we recover the product from [23]

$$\frac{\pi}{2} = \left(\frac{2}{1}\right)^{1/2} \left(\frac{2^2}{1 \cdot 3}\right)^{1/4} \left(\frac{2^3 \cdot 4}{1 \cdot 3^3}\right)^{1/8} \left(\frac{2^4 \cdot 4^4}{1 \cdot 3^6 \cdot 5}\right)^{1/16} \cdots.$$

**Example 5.2.** Setting $m = 1$, $z = -1$ and $u = 1$ and multiplying by 4, we get

$$4\frac{\partial \Phi}{\partial s}(-1, -1, 1) = \sum_{n=0}^{\infty} \frac{n+1}{2^n} \sum_{k=0}^n (-1)^{k+1} \binom{n}{k} \ln(k+1).$$

Using (18), we obtain the product
$$\frac{A^{12}}{2^{4/3}e} = \left(\frac{2}{1}\right)^{2/2} \left(\frac{2^2}{1 \cdot 3}\right)^{3/4} \left(\frac{2^3 \cdot 4}{1 \cdot 3^3}\right)^{4/8} \left(\frac{2^4 \cdot 4^4}{1 \cdot 3^6 \cdot 5}\right)^{5/16} \cdots,$$
which converges faster than the products involving $A$ in Examples 5.7 and 5.11.

**Example 5.3.** Setting $m = 2$, $z = -1$ and $u = 1$, we get
$$\frac{\partial \Phi}{\partial s}(-1, -2, 1) = \sum_{n=0}^{\infty} \frac{n^2 + n}{2^{n+3}} \sum_{k=0}^{n} (-1)^{k+1} \binom{n}{k} \ln(k+1).$$

Then (19) gives
$$e^{7\zeta(3)/4\pi^2} = \left(\frac{2}{1}\right)^{1/8} \left(\frac{2^2}{1 \cdot 3}\right)^{3/16} \left(\frac{2^3 \cdot 4}{1 \cdot 3^3}\right)^{6/32} \left(\frac{2^4 \cdot 4^4}{1 \cdot 3^6 \cdot 5}\right)^{10/64} \cdots.$$

**Example 5.4.** Take $m = 0$, $z = -1$ and $u = 1/2$. Writing $\ln\left(k + \frac{1}{2}\right)$ as $\ln(2k+1) - \ln 2$ and using (25), we see that
$$\frac{\partial \Phi}{\partial s}\left(-1, 0, \frac{1}{2}\right) - \frac{\ln 2}{2} = \sum_{n=0}^{\infty} \frac{1}{2^{n+1}} \sum_{k=0}^{n} (-1)^{k+1} \binom{n}{k} \ln(2k+1).$$

Using (28), we get the product
$$\frac{\Gamma(1/4)}{2\Gamma(3/4)} = \left(\frac{3}{1}\right)^{1/4} \left(\frac{3^2}{1 \cdot 5}\right)^{1/8} \left(\frac{3^3 \cdot 7}{1 \cdot 5^3}\right)^{1/16} \left(\frac{3^4 \cdot 7^4}{1 \cdot 5^6 \cdot 9}\right)^{1/32} \cdots.$$

**Example 5.5.** Setting $m = 1$, $z = -1$ and $u = 1/2$, we get
$$\frac{\partial \Phi}{\partial s}\left(-1, -1, \frac{1}{2}\right) = \sum_{n=0}^{\infty} \frac{n}{2^{n+2}} \sum_{k=0}^{n} (-1)^{k+1} \binom{n}{k} \ln(2k+1).$$

Use (20) and exponentiate to obtain
$$e^{G/\pi} = \left(\frac{3}{1}\right)^{1/8} \left(\frac{3^2}{1 \cdot 5}\right)^{2/16} \left(\frac{3^3 \cdot 7}{1 \cdot 5^3}\right)^{3/32} \left(\frac{3^4 \cdot 7^4}{1 \cdot 5^6 \cdot 9}\right)^{4/64} \cdots.$$

**Corollary 5.2.** *For all complex $s \neq 1$,*
$$\Phi(1, s, u) + (s-1)\frac{\partial \Phi}{\partial s}(1, s, u) = \sum_{n=0}^{\infty} \frac{1}{n+1} \sum_{k=0}^{n} (-1)^{k+1} \binom{n}{k} \frac{\ln(u+k)}{(u+k)^{s-1}}.$$

*Proof.* Multiply (23) by $s - 1$ and differentiate with respect to $s$. □

**Example 5.6.** Set $u = 1$ and $s = 0$ and multiply by 2. The relation (2) and the values $\zeta(0) = -1/2$ and $\zeta'(0) = -\frac{1}{2}\ln 2\pi$ give
$$\ln 2\pi - 1 = \sum_{n=0}^{\infty} \frac{2}{n+1} \sum_{k=0}^{n} (-1)^{k+1} \binom{n}{k} (k+1) \ln(k+1).$$

Exponentiating, we obtain the product
$$\frac{2\pi}{e} = \left(\frac{2^2}{1}\right)^{2/2} \left(\frac{2^4}{1 \cdot 3^3}\right)^{2/3} \left(\frac{2^6 \cdot 4^4}{1 \cdot 3^9}\right)^{2/4} \left(\frac{2^8 \cdot 4^{16}}{1 \cdot 3^{18} \cdot 5^5}\right)^{2/5} \cdots.$$

**Example 5.7.** Set $s = -1$ and $u = 1$ and multiply by $1/2$. Using relations (2) and (17), and the value $\zeta(-1) = -1/12$, we get the formula

$$\ln A - \frac{1}{8} = \sum_{n=0}^{\infty} \frac{1}{2n+2} \sum_{k=0}^{n} (-1)^{k+1} \binom{n}{k} (k+1)^2 \ln(k+1), \tag{54}$$

which yields

$$\frac{A}{e^{1/8}} = \left(\frac{2^4}{1}\right)^{1/4} \left(\frac{2^8}{1 \cdot 3^9}\right)^{1/6} \left(\frac{2^{12} \cdot 4^{16}}{1 \cdot 3^{27}}\right)^{1/8} \cdots.$$

**Theorem 5.1.** *For $u > 0$,*

$$\psi(u) = \sum_{n=0}^{\infty} \frac{1}{n+1} \sum_{k=0}^{n} (-1)^k \binom{n}{k} \ln(u+k). \tag{55}$$

*Proof.* Multiply (23) by $s-1$ and compute the derivative with respect to $s$ at $s = 1$, using (42). □

**Example 5.8.** Take $u = 1$, multiply by $-1$, and exponentiate. Using (53), we recover Ser's product [17] (rediscovered in [21] and [23])

$$e^\gamma = \left(\frac{2}{1}\right)^{1/2} \left(\frac{2^2}{1 \cdot 3}\right)^{1/3} \left(\frac{2^3 \cdot 4}{1 \cdot 3^3}\right)^{1/4} \left(\frac{2^4 \cdot 4^4}{1 \cdot 3^6 \cdot 5}\right)^{1/5} \cdots.$$

Compare the remarkably similar product for $e$ in Example 5.12.

**Example 5.9.** Take first $u = 1/4$ and then $u = 3/4$. Using the relation $\psi(3/4) - \psi(1/4) = \pi$ (see [24, p. 427]), we obtain the series

$$\pi = \sum_{n=0}^{\infty} \frac{1}{n+1} \sum_{k=0}^{n} (-1)^k \binom{n}{k} \ln \frac{4k+3}{4k+1},$$

which gives the product

$$e^\pi = \left(\frac{3}{1}\right)^{\frac{1}{1}} \left(\frac{3 \cdot 5}{1 \cdot 7}\right)^{\frac{1}{2}} \left(\frac{3 \cdot 5^2 \cdot 11}{1 \cdot 7^2 \cdot 9}\right)^{\frac{1}{3}} \left(\frac{3 \cdot 5^3 \cdot 11^3 \cdot 13}{1 \cdot 7^3 \cdot 9^3 \cdot 15}\right)^{\frac{1}{4}} \left(\frac{3 \cdot 5^4 \cdot 11^6 \cdot 13^4 \cdot 19}{1 \cdot 7^4 \cdot 9^6 \cdot 15^4 \cdot 17}\right)^{\frac{1}{5}} \cdots.$$

For the next result, recall the formulas for the *Euler beta function* [25, p. 254],

$$B(u,v) = \int_0^1 x^{u-1}(1-x)^{v-1} dx = \frac{\Gamma(u)\Gamma(v)}{\Gamma(u+v)}.$$

**Theorem 5.2.** *For $j = 1, 2, \ldots$,*

$$B(u,j) = \sum_{n=j}^{\infty} \frac{1}{n-j+1} \sum_{k=0}^{n} (-1)^{k+1} \binom{n}{k} \ln(u+k). \tag{56}$$

*Proof.* We first establish the formula

$$\int_0^1 \frac{x^{u-1}(1-x)^n}{-\ln x} dx = \sum_{k=0}^{n} (-1)^{k+1} \binom{n}{k} \ln(u+k) \tag{57}$$

for $n = 1, 2, \ldots$ . Fix $n$ and define $f(u)$ to be the difference between the integral and the sum in (57). Then $f'(u) = 0$. (To see this, differentiate $f(u)$ under the integral sign, expand $(1-x)^n$, and evaluate the integral.) Thus $f(u) = C$ is a constant. By (25) with $m = 0$, we may replace $\ln(u+k)$ with $\ln\left(1 + \frac{k}{u}\right)$ in $f(u)$. It follows that $f(u) \to 0$ as $u \to \infty$. Hence $C = 0$ and (57) follows. Now to prove (56), multiply (57) by $(n-j+1)^{-1}$ and sum from $n = j$ to $\infty$. $\square$

**Example 5.10.** Set $u = 1$ in (57), multiply the equation by $(-1)^{n+1}$, and sum from $n = 1$ to $\infty$. Using Example 3.13 (let $x = X/Y$, $y = Y$ and integrate with respect to $Y$), we obtain

$$\ln \sigma = \int_0^1 \frac{1-x}{(2-x)(-\ln x)} dx = \sum_{n=1}^{\infty} (-1)^n \sum_{k=0}^{n} (-1)^k \binom{n}{k} \ln(k+1),$$

where $\sigma$ is Somos's constant (37). Exponentiation gives a product for $\sigma$ with rational factors:

(58) $$\sigma = \frac{2}{1} \cdot \frac{1 \cdot 3}{2^2} \cdot \frac{2^3 \cdot 4}{1 \cdot 3^3} \cdot \frac{1 \cdot 3^6 \cdot 5}{2^4 \cdot 4^4} \cdots .$$

Since the partial products are alternately above and below $\sigma$, and converge slowly, the product is a good candidate for acceleration by Euler's transformation [19], [23]. Applying it and simplifying (as in [23, section 3]), we arrive at the product

(59) $$\sigma = \left(\frac{2}{1}\right)^{1/2} \left(\frac{3}{2}\right)^{1/4} \left(\frac{4}{3}\right)^{1/8} \left(\frac{5}{4}\right)^{1/16} \cdots ,$$

which converges faster than the products for $\sigma$ in both (37) and (58). (Note that (59) also follows by writing $\sigma$ as $\sigma^2/\sigma$ and substituting the product in (37).)

**Example 5.11.** Set $u = 1$ and $j = 4$ in (56), and replace $n$ with $n+4$. Multiplying the equation by $1/2$, so that the left side becomes $B(1,4)/2 = 1/8$, we add the result to (54). Replacing $k$ with $k-1$, and $n$ with $n-1$, we obtain

$$\ln A = \sum_{n=1}^{\infty} \frac{1}{2n} \left[ \sum_{k=1}^{n} (-1)^k \binom{n-1}{k-1} k^2 \ln k + \sum_{k=1}^{n+4} (-1)^k \binom{n+3}{k-1} \ln k \right].$$

This yields a product for the Glaisher-Kinkelin constant:

$$A = \left(\frac{2^4 \cdot 4^4}{1 \cdot 3^6 \cdot 5}\right)^{1/2} \left(\frac{2^9 \cdot 4^{10} \cdot 6}{1 \cdot 3^{10} \cdot 5^5}\right)^{1/4} \left(\frac{2^{14} \cdot 4^{20} \cdot 6^6}{1 \cdot 3^{24} \cdot 5^{15} \cdot 7}\right)^{1/6} \cdots .$$

The remaining results use only the series for the Euler beta function in Theorem 5.2. In particular, their construction does not involve Lerch's function.

**Example 5.12.** As $B(1,j) = 1/j$, exponentiating (56) with $u = 1$ yields a product for $e^{1/j}$. When $j = 1$, this is the first author's product for $e$ in [23]:

$$e = \left(\frac{2}{1}\right)^{1/1} \left(\frac{2^2}{1 \cdot 3}\right)^{1/2} \left(\frac{2^3 \cdot 4}{1 \cdot 3^3}\right)^{1/3} \left(\frac{2^4 \cdot 4^4}{1 \cdot 3^6 \cdot 5}\right)^{1/4} \cdots .$$

The extension to $e^{1/j}$ for $j = 2, 3, \ldots$ was first found by C. Goldschmidt and J. B. Martin [8] using probability theory. With $j = 2$, it gives

$$\sqrt{e} = \left(\frac{2^2}{1 \cdot 3}\right)^{1/1} \left(\frac{2^3 \cdot 4}{1 \cdot 3^3}\right)^{1/2} \left(\frac{2^4 \cdot 4^4}{1 \cdot 3^6 \cdot 5}\right)^{1/3} \cdots .$$

**Theorem 5.3.** *The infinite product representation of the exponential function*

$$(60) \qquad e^x = \prod_{n=1}^{\infty} \left( \prod_{k=1}^{n} (kx+1)^{(-1)^{k+1} \binom{n}{k}} \right)^{1/n}$$

*is valid for any real number $x \geq 0$.*

*Proof.* The formula holds when $x = 0$. If $x > 0$, then since $u$, unlike $j$, is not restricted to integer values in Theorem 5.2, we may take $u = 1/x$. By (25) with $m = 0$, we may then replace $\ln\left(\frac{1}{x} + k\right)$ with $\ln(kx+1)$ in (56). Setting $j = 1$ and using $B(1/x, 1) = x$, we exponentiate (56) and obtain (60). □

**Example 5.13.** Taking $x = 2$, we get

$$e^2 = \left(\frac{3}{1}\right)^{1/1} \left(\frac{3^2}{1 \cdot 5}\right)^{1/2} \left(\frac{3^3 \cdot 7}{1 \cdot 5^3}\right)^{1/3} \left(\frac{3^4 \cdot 7^4}{1 \cdot 5^6 \cdot 9}\right)^{1/4} \cdots.$$

When $x = p/q$ is a rational number with $p > 0$ and $q > 0$, we can simplify (60) by replacing $kx + 1$ with $kp + q$, again using (25) with $m = 0$. An example is

$$e^{2/3} = \left(\frac{5}{3}\right)^{1/1} \left(\frac{5^2}{3 \cdot 7}\right)^{1/2} \left(\frac{5^3 \cdot 9}{3 \cdot 7^3}\right)^{1/3} \left(\frac{5^4 \cdot 9^4}{3 \cdot 7^6 \cdot 11}\right)^{1/4} \cdots.$$

In particular, with $x = 1/q$ and $q = 2, 3, \ldots$, the roots $e^{1/q}$ are represented as different products from those in Example 5.12. For instance,

$$\sqrt{e} = \left(\frac{3}{2}\right)^{1/1} \left(\frac{3^2}{2 \cdot 4}\right)^{1/2} \left(\frac{3^3 \cdot 5}{2 \cdot 4^3}\right)^{1/3} \left(\frac{3^4 \cdot 5^4}{2 \cdot 4^6 \cdot 6}\right)^{1/4} \cdots.$$

**Acknowledgments.** For bringing Ser's and Hasse's papers to our attention, we thank Donal Connon and Alex Schuster, respectively.

Av. Cesáreo Alierta 31, esc. izda 4ºA
Zaragoza 50008 Spain
*E-mail address*: jguillera@able.es
*URL*: http://personal.auna.com/jguillera

209 West 97th Street Apt 6F
New York, NY 10025 USA
*E-mail address*: jsondow@alumni.princeton.edu
*URL*: http://home.earthlink.net/j̃sondow/